\hfuzz=6pt
\font\titlefont=cmbx10 scaled\magstep1

\magnification=1200
\line{}
\vskip 1.5cm
\centerline{\titlefont AN ALGEBRAIC INTERPRETATION OF THE}
\smallskip
\centerline{\titlefont CONTINUOUS BIG q-HERMITE POLYNOMIALS}
\vskip 2cm
\centerline{\bf Roberto Floreanini}
\smallskip
\centerline{Istituto Nazionale di Fisica Nucleare, Sezione di Trieste}
\centerline {Dipartimento di Fisica Teorica,
Universit\`a di Trieste}
\centerline{Strada Costiera 11, 34014 Trieste, Italy}
\vskip 1cm
\centerline{\phantom{$^{(*)}$}{\bf Jean LeTourneux}$^{(*)}$}
\medskip
\centerline{and}
\medskip
\centerline{\phantom{$^{(*)}$}{\bf Luc Vinet}\footnote{$^{(*)}$}{Supported 
in part by the National Sciences and Engineering
Research Council \hbox{(NSERC)} of Canada and the Fonds FCAR of Qu\'ebec.}}
\smallskip
\centerline{Centre de Recherches Math\'ematiques}
\centerline{Universit\'e de Montr\'eal}
\centerline{Montr\'eal, Canada H3C 3J7}
\vskip 2cm
\centerline{\bf Abstract}
\smallskip\midinsert\narrower\narrower\noindent  
The continuous big $q$-Hermite polynomials are shown to realize a basis
for a representation space of an extended $q$-oscillator algebra.
An expansion formula is algebraically derived using this model.
\endinsert

\vfill\eject

Lie algebra theory is well known to provide a unifying framework for discussing
special functions.  The discovery, some ten years ago, of quantum groups has in
turn prompted the undertaking of a systematic investigation of the algebraic
properties of the 
$q$-analogs of those special functions.  One indeed witnesses nowadays intense
research activity in this area as $q$-special functions are  seen to have more
and more applications.

Within the Askey scheme,$^{1,2}$ sets of basic or $q$-orthogonal polynomials
are called continuous because their elements are orthogonal with respect to
continuous measures.  We have initiated in Ref.~[3] a study of these continuous
$q$-polynomials from an algebraic point of view, 
focusing on the continuous $q$-Hermite and
continuous $q$-ultraspherical polynomials and, as a result, have shown that
various properties of these functions can be derived using symmetry techniques. 
We indicate here that the class of continuous big $q$-Hermite
polynomials also lends itself to a similar treatment.

We shall be using standard notation.$^{1,2}$ The $q$-hypergeometric series
${}_r \phi_s$ is
$$
{}_r \phi_s \bigg({a_1, a_2, \dots , a_r\atop b_1, \dots , b_s}\bigg| q; z 
\bigg)
= \sum_{n=0}^\infty {(a_1,\dots , a_r; q)_n\over (q,b_1,\dots ,b_s;q)_n} 
\bigl[(-1)^n q^{n(n-1)/2}\bigr]^{1+s-r} z^n\ , \eqno(1)
$$
with
$$
(a_1,a_2, \dots , a_k; q)_\alpha = (a_1;q)_\alpha (a_2;q)_\alpha \dots
(a_k;q)_\alpha\ ,\eqno(2a)
$$
and
$$
(a;q)_\alpha ={(a;q)_\infty\over (aq^\alpha;q)_\infty}\ ,\qquad
(a;q)_\infty = \prod_{k=0}^\infty (1-aq^k)\ ,  \quad |q| < 1\ . \eqno(2b)
$$
Clearly, the series ${}_r\phi_s$ terminates if one of the $a_i$, $i=1, \dots ,
r$, is equal to $q^{-n}$ with $n$ a positive integer.

The continuous big $q$-Hermite polynomials  $H_n(x; a\vert  q)$ depend on one
parameter and are defined as follows:$^2$
$$
\eqalignno{
H_n(x; a\vert  q) 
	&= a^{-n} {}_3\phi_2\bigg({q^{-n},\atop\ }{ae^{i\theta},\atop 0,}
{ae^{-i\theta}\atop 0}\bigg| q;q \bigg) &(3a)\cr
 	& = e^{in\theta} {}_2\phi_0
\bigg({q^{-n}, ae^{i\theta}\atop -}\bigg| q;q^n e^{-2i\theta}\bigg)\ ,
\qquad x = \cos \theta\ . &(3b)}
$$
When $a$ is real and $|a| < 1$, these polynomials obey the following
orthogonality relation:
$$
{1\over 2\pi} \int_{-1}^1 {w(x;a\vert  q)\over \sqrt{1-x^2}} H_m(x;a\vert  q)
H_n(x;a\vert  q)\, dx = {\delta_{mn}\over (q^{n+1}; q)_\infty}\ , \eqno(4)
$$
where
$$
w(x;a\vert  q) =
\left|{(e^{2i\theta};q)_\infty\over (ae^{i\theta};q)_\infty} 
\right|^2\ . \eqno(5)
$$
Note that as $q \rightarrow 1^-$,
$$
H_n(x;a\vert  q) \rightarrow (2x - a)^n\ .\eqno(6)
$$
The continuous $q$-Hermite polynomials  $H_n(x\vert  q)$ 
can be defined as the $a\rightarrow 0$ 
limit of the polynomials $H_n(x;a\vert  q)$.  This limit can be taken
immediately in $(3b)$, leading to
$$
H_n(x\vert  q) =
	e^{in\theta}\ {}_2\phi_0
\bigg({q^{-n},\ 0\atop -}\bigg| q; q^n e^{-2i\theta}\bigg)\ . \eqno(7)
$$
The orthogonality relations of these polynomials are obtained by 
setting $a=0$ in
(4) and (5) and writing $H_n(x;0\vert  q) \equiv H_n(x\vert  q)$.

In the following, we shall show that the continuous big $q$-Hermite polynomials
occur in the realization of a 
set of basis vectors for a representation space of
a $q$-algebra ${\cal G}_q$ that encompasses the $q$-Heisenberg algebra.  This
algebraic set up will then be used to derive an 
expansion formula involving the polynomials $H_n(x;a\vert  q)$.

We shall now present a realization of ${\cal G}_q$ in terms of operators acting on
functions of the two variables: $x = (z + z^{-1})/2$, 
with $z = e^{i\theta}$, and
$t$.  To do so, we shall need the $q$-shift operators $T_z$ and $T_t$ whose
powers act as follows:
$$
\eqalign{
&T_z^\alpha f \bigl[ (z+z^{-1}), t \bigr]
	= f \bigl[ (q^\alpha z + q^{-\alpha} z^{-1}), t \bigr]\ ,\cr
&T_t^\beta f \bigl[ (z+z^{-1}), t \bigr]
	= f \bigl[ (z + z^{-1}), q^{\beta}t  \bigr]\ , \qquad\quad \alpha , 
\beta \in {\bf R}\ .}\eqno(8)
$$
Let,
$$
\eqalignno{
&A_+
	= {t\over z-z^{-1}} \big(T_z^{1/2} - T_z^{-1/2}\big)\ , &(9a)\cr
&A_-
	= {q^{-1/2}\over t(z-z^{-1})} \biggl[
{1\over z^2} (1-q^{-1/2} z \, T_t^{1/2}) \, T_z^{1/2} \cr
& \hskip 6cm -z^2 \left(1 -{q^{-1/2}\over z} T_t^{1/2}
\right) T_z^{-1/2} \biggr]\ , &(9b)\cr
&B_+
	= {t\over (z-z^{-1})} (z \, T_z^{-1/2} - {1\over z} \, T_z^{1/2})\ ,
&(9c)\cr
&B_-
	= {1\over t(z-z^{-1})} 
 \biggl[ z
\left(1-{q^{-1/2}\over z} T_t^{1/2} \right) T_z^{-1/2} \cr
&\hskip 6cm -{1\over z} \left( 1 - q^{-1/2} z \,
T_t^{1/2}\right) \,  T_z^{1/2} \biggr]\ , &(9d)\cr
&K = T_t\ . &(9e)}
$$
Notice that as $q \rightarrow 1^-$,
$$
{1\over1-q} A_+ \rightarrow -{t\over 2} {\partial\over\partial x}, 
\quad  A_- \rightarrow -{1\over t} (2x-1), \quad B_+ \rightarrow t,
\quad B_- \rightarrow {1\over t}\ , \eqno(10a)
$$
and
$$
{1-K\over 1-q} \rightarrow t {\partial\over \partial t}\ . \eqno(10b)
$$
In this limit, $B_+$, $B_-$ and $K$ enlarge in a simple way the Heisenberg
algebra that $A_+/(1-q)$ and $A_-$ realize.  Consider the set of functions
$$
f_n^m(x,t) = t^mH_n (x; q^{m/2} \vert  q)\ , \qquad n \in {\bf Z}^+\ , 
\qquad m \in {\bf Z}\ . \eqno(11)
$$
It can be checked$^4$ that the operators $A_+$, $A_-$, $B_+$, $B_-$ and $K$
transform this ensemble of functions onto itself according to:
$$
\eqalignno{
&A_+\, f_n^m = -q^{-n/2} (1-q^n) f_{n-1}^{m+1}\ , &(12a)\cr
&A_-\, f_n^m= -q^{-(n+1)/2} f_{n+1}^{m-1}\ , &(12b)\cr
&B_+\,  f_n^m= q^{-n/2} f_n^{m+1}\ , &(12c)\cr
&B_-\, f_n^m= q^{-n/2} f_n^{m-1}\ , &(12d)\cr
&K\, f_n^m= q^m f_n^m\ . &(12e)}
$$
It is also natural to consider two additional operators, namely, multiplication
by $x$ and by $t^2$.  Take $P = 2x$ and $Q = t^2$.  The three-term recurrence
relation of the continuous big $q$-Hermite polynomials$^2$
$$
2x H_n(x;a \vert  q) =
	H_{n+1}(x; a \vert  q) + 
a q^n H_n (x;a \vert  q) + (1 - q^n) H_{n-1}(x;a\vert  q)\ ,\eqno(13)
$$
immediately gives the action of $P=2x$ on the space of functions spanned by the
$f_n^m$.  It reads
$$
Pf_n^m =
	f_{n+1}^m + q^{n+m/2} f_n^m + (1 - q^n) f_{n-1}^m\ .\eqno(14)
$$
In order to write down the action of $Q = t^2$ on the basis functions, one first
observes that the continuous big $q$-Hermite polynomials satisfy the following
identity
$$
H_n (x;a \vert  q) =
	H_n(x;aq \vert  q) - a(1 - q^n) \,  H_{n-1}(x;aq \vert  q)\ .\eqno(15)
$$

This formula is most easily proven by checking that both sides verify the same
recurrence relation with the same initial condition.  It then follows from (15)
that
$$
Qf_n^m =
	f_n^{m+2} - q^{m/2} (1-q^n) f_{n-1}^{m+2}\ .\eqno(16)
$$
The $q$-algebra that the operators $A_\pm$, $B_\pm$, $K$, $P$ and $Q$ 
realize can
be characterized by the following relations:
$$
\eqalign{
&A_- A_+ - qA_+A_- = -(1-q)\ ,\cr
&B_+A_+ - q^{1/2} A_+B_+=0\ ,\cr
&A_-B_+ -q^{1/2} B_+A_- =0\ ,\cr
&A_+P - q^{1/2}PA_+ = -q^{-1/2}(1-q)B_+\ ,\cr
&q^{1/2}B_+P - PB_+= (1-q)A_+\ ,\cr
&A_+Q - QA_+= 0\ , \cr
&B_+Q - QB_+=0\ , \cr
&KA_+ - qA_+K=0\ , \cr
&KB_+ - qB_+K=0\ ,\cr
&KP - PK=0\ ,}\qquad
\eqalign{
&B_+B_- - B_-B_+ = 0\ ,\cr
&B_-A_+ - q^{1/2} A_+B_-= 0\ ,\cr
&A_-B_- - q^{1/2}B_-A_- = 0\ ,\cr
&A_-P - q^{-1/2} PA_-= q^{-1}(1-q)B_-\ ,\cr
&B_-P - q^{1/2}PB_-= -(1-q)A_-\ ,\cr
&(QA_- - A_-Q) =q^{-1}(1-q)B_+K^{1/2}\ ,\cr
&B_-Q - q \, QB_-= (1-q)B_+\ ,\cr
&KA_- - q^{-1}A_- K=0\ ,\cr
&KB_- - q^{-1} B_- K= 0\ ,\cr
&KQ - q^2 QK= 0\ .}\eqno(17) 
$$

Let us make a few comments on this algebra.  First note that $A_+$ and $A_-$
generate the $q$-Heisenberg algebra.$^5$  There are also various 
interesting $q$-subalgebras.  
The generators $A_+$, $B_+$ and $P$, for example, form a 
closed set.  We see that $A_+$ and $B_+$ $q$-commute and 
are in a certain way rotated one into the
other by $P$.  The set $\{ A_-, B_-, P\}$ also has a similar structure.  
When $q\rightarrow 1$, the algebra exhibits large abelian sectors. 
Furthermore, some generators (see (10)) become redundant: 
$B_+^2$ and $Q$ for instance, 
have the same limit, and the same is
true of $B_+A_-$ and $-2P + 1$.

We now want to illustrate how the model given in (9) and (11)
can be used to derive properties of the 
continuous big $q$-Hermite polynomials.  In the Lie theory
approach to ordinary special functions, one considers exponentials of 
the algebra generators and 
relates their matrix elements in representation spaces to various functions of
interest.  One then uses diverse realizations to obtain identities 
and formulas. 
To proceed similarly in the case of $q$-special functions, we need 
$q$-analogs  of the exponential.  
It has been appreciated that the $q$-exponential which is naturally
associated to the continuous $q$-orthogonal polynomials is the one first
introduced in Ref.[6] and denoted by ${\cal E}_q(x;a,b)$.  Indeed, one for
example finds$^{6,7,3}$ that it generates the 
continuous $q$-Hermite polynomials:
$$
{\cal E}_q \left(x; -i, b/2 \right) =
	\left(- b^2/4; q^2 \right)_\infty^{-1} \sum_{k=0}^\infty
{q^{k^2/4}\over (q;q)_k}\, \left({ib\over 2}\right)^k H_k(x \vert  q)\ .
\eqno(18)
$$
We shall analogously consider the ${\cal E}_q$-exponential 
of the generator $P=2x$
and determine some of its matrix elements in the bases $\{ f_n^m \}$ 
to obtain an interesting expansion 
formula in the polynomials $H_n(x; q^m \vert  q)$.  However before doing 
so, we need to define this $q$-exponential ${\cal E}_q$ and to record some of
its properties.

We have$^6$
$$
{\cal E}_q(x;a,b) =
	\sum_{n=0}^\infty {q^{n^2/4}\over(q;q)_n} \,
\bigl( aq^{(1-n)/2} e^{i\theta};
q \bigr)_n \bigl( aq^{(1-n)/2} e^{-i\theta}; q\bigr)_n \,b^n, \qquad x =
\cos\theta\ . \eqno(19)
$$
In the limit $q \rightarrow 1^-$,
$$
{\cal E}_q \bigl(x;a,(1-q)b\bigr) \rightarrow 
\exp\bigl[ (1+a^2 - 2a\,x)b \bigr]\ ,\eqno(20)
$$
and, in particular, for $a=-i$,
$$
\lim_{q \rightarrow 1^-} {\cal E}_q \left( x; -i, (1-q) b/2 \right) = e^{ibx}
\ . \eqno(21)
$$
The essential feature of these $q$-exponentials is that they are eigenfunctions
of the divided difference operator
$$
\tau = {1\over z-z^{-1}} (T_z^{1/2} - T_z^{-1/2})\ . \eqno(22)
$$
Indeed,
$$
\tau \,{\cal E}_q(x;a,b) = ab\, q^{-1/4}\, {\cal E}_q(x;a,b)\ . \eqno(23)
$$
The continuous $q$-polynomials obey second order $\tau$-difference equations. 
Note also that in our model $A_+ = t \tau$ (see $(9a)$).  It is thus not
surprising, in view of (23), to see ${\cal E}_q$ be the appropriate 
$q$-exponential to use in connection 
with continuous $q$-polynomials.  
There is one more property of ${\cal E}_q$ that
we shall need in the following.  Consider the function $g_n(b)$ defined by
$$
g_n(b) = {\cal E}_q \bigl(-;0,b\,q^{n/2} \bigr) = \sum_{k=0}^\infty
{q^{k(k+2n)/4}\over(q;q)_k}\, b^k\ . \eqno(24)
$$
It is readily verified that $g_n(b)$ satisfies the $3$-term recurrence relation
$$
g_{n+1}(b) = g_{n-1} (b) - b\, q^{(2n-1)/4}\, g_n(b)\ .\eqno(25)
$$

As an example of application of our formalism, we will now derive an expansion
formula for ${\cal E}_q(x;-i,b/2)$ in terms of continuous big $q$-Hermite
polynomials.  This $q$-exponential of the generator $P/2 = x$ acts, 
of course, on the representation space 
of our $q$-algebra.  Recall that $K = T_t$ is diagonal on the basis $\{ f_n^m
\}$: $K\, f_n^m = q^m f_n^m$.  Since $P$ and $K$ commute, we must have
$$
{\cal E}_q \left( x; -i, b/2 \right)\, f_0^m = \sum_{n=0}^\infty 
W_n^m (b)\, f_n^m\ .\eqno(26)
$$
Note that we are considering the action of the $q$-exponential of $x$ on the
particular basis vectors $f_0^m(x,t) = t^m$.  The expansion coefficients
$W_n^m(b)$ will be obtained from the recursion relations that they obey.  These
relations will be found by 
exploiting properties of the ${\cal E}_q$-exponential and making use of the
representation (12), (15), (16).  Let us first act on both sides of (26) with
$A_+ = t \tau$.  With the help of (23), we see on the one hand that
$$
\eqalign{
A_+ {\cal E}_q \left( x; -i, b/2 \right) f_0^m
	&= -i {b\over 2}\, q^{-1/4} {\cal E}_q \left( x; -i, b/2 \right) 
f_0^{m+1}\cr
	&= -i{b\over 2}\, q^{-1/4} \sum_{n=0}^\infty W_n^{m+1}(b)\, f_n^{m+1}\ ,}
\eqno(27)
$$
the last equality following from (26).  On the other hand, using $(12a)$, 
we have
$$
\eqalign{
A_+ {\cal E}_q \left( x; -i, b/2 \right)\, f_0^m
	&=\sum_{n=0}^\infty W_n^m (b) A_+\, f_n^m\cr
	&= -\sum_{n=0}^\infty q^{-n/2} (1-q^n) W_n^m(b)\, f_{n-1}^{m+1}\ .} 
\eqno(28)
$$
Equating the right-hand sides of (27) and (28), we then find 
$$
i (b/2) q^{-1/4}\,  W_n^{m+1} (b)
	= q^{-(n+1)/2} (1-q^{n+1})\,  W_{n+1}^m (b)\ .\eqno(29)
$$
Second, we act similarly on both sides of (26) with $Q = t^2$. Clearly,
$$
\eqalign{
Q \, {\cal E}_q \left( x; -i, b/2 \right) 
	&= {\cal E}_q \left( x; -i, b/2 \right) f_0^{m+2}\cr
	&= \sum_{n=0}^\infty W_n^{m+2}(b) \, f_n^{m+2}\ ,}
\eqno(30)
$$
while (16) yields
$$
\eqalign{
Q\, {\cal E}_q \left( x; -i, b/2 \right)\, f_0^m 
	&= \sum_{n=0}^\infty W_n^{m}(b)\,  Q\, f_n^{m}\cr
	&= \sum_{n=0}^\infty W_n^{m}(b)\, \bigl[ f_n^{m+2} - q^{m/2}
(1-q^n)f_{n-1}^{m+2} \bigr]\ .}
\eqno(31)
$$
Combining (30) and (31), we get
$$
W_n^{m+2}(b)
	= W_n^m(b) - (1 - q^{n+1}) q^{m/2}\, W_{n+1}^m(b)\ . \eqno(32)
$$
Finally, we replace in this last equation $W_{n+1}^m (b)$ by the expression 
that (29) gives for it to find
$$
W_n^{m+2}(b) 
= W_n^m (b) - i {b\over 2}\, q^{(n+m+1/2)/2}\, W_n^{m+1}(b)\ . \eqno(33)
$$
The matrix elements $W_n^m(b)$ can now be explicitly determined from the two
recurrence relations (29) and (33), that we have found for them.  Separation of
the discrete variables is readily achieved in these equations by taking
$W_n^m(b)$ of the form
$$
W_n^m(b)= u_n(b)\, y_{m+n}(b)\ . \eqno(34)
$$
Indeed, substitution of (34) in (29) and (33), respectively, gives
$$
u_{n+1} (b) 
= {q^{(2n+1)/4}\over 1-q^{n+1}}\,  \left({ib\over2}\right)\, u_n(b)\ ,\eqno(35)
$$
and
$$
y_{m+n+2} (b) 
	= y_{m+n}(b) -{ib\over 2}\, q^{(m+n+1/2)/2}\, y_{m+n+1}(b)\ . \eqno(36)
$$
The recurrence relation (35) is easily solved and fixes $u_n(b)$ up to a
function $u_0(b)$:
$$
u_n(b) 
= u_0(b)\, {q^{n^2/4}\over(q;q)_n} \left({ib\over2}\right)^n\ . \eqno(37)
$$
The $3$-term recurrence relation (36) is recognized as the one already given in
(25) and $y_{m+n}(b)$ is thus immediately identified:
$$
y_{m+n}(b) 
= {\cal E}_q \left(-; 0, {ib\over2} q^{(m+n)/2} \right)\ . \eqno(38)
$$
(It is understood that the overall arbitrary function of $b$ in the solution of
(36), is to be absorbed in $u_0(b)$.)  If we use the realization $f_n^m(x,t) =
t^m H_n \bigl( x;q^{m/2} \vert  q \bigr)$ and factor out 
the $t$-dependence, (26) becomes
$$
{\cal E}_q \left( x; -i, b/2 \right)
	= \sum_{n=0}^\infty W_n^m(b)\, H_n\bigl( x; q^{m/2} \vert  q \bigr)\ .
\eqno(39)
$$
At this point,
$$
W_n^m(b)
= u_0(b)\, {q^{n^2/4}\over(q;q)_n}\, {\cal E}_q \left( -; 0,{ib\over2}
q^{(m+n)/2} \right)\, \left({ib\over2}\right)^n\ .\eqno(40)
$$
We will therefore have obtained the identity we are looking for, once we will
have determined $u_0(b)$.  
To this end, notice that as $m \rightarrow \infty$, or $q^m\rightarrow 0$,
$$
W_n^m(b) \rightarrow u_0(b)\, {q^{n^2/4}\over(q;q)_n}
\left({ib\over2}\right)^n\ ; \eqno(41)
$$
recall also that $H_n(x; 0 \vert  q) \equiv H_n(x \vert  q)$.  Hence, in the
limit $q^m \rightarrow 0$, (39) must coincide with the expansion formula for 
${\cal E}_q(x; -i, b/2)$ 
in terms of continuous $q$-Hermite polynomials already given in (18). 
Comparison immediately shows that
$$
u_0(b) = (-b^2/4; q^2)^{-1}_\infty\ .\eqno(42)
$$
Putting everything together finally gives the 
following expansion formula of the
${\cal E}_q$-expo-\break 
nential of $x$ in terms of continuous big $q$-Hermite 
polynomials:
$$
\eqalign{
{\cal E}_q (x; -i, b/2) =(-b^{2/4}; q^2)^{-1}_\infty\, 
\sum_{n=0}^\infty {q^{n^2/4}\over(q;q)_n}\,
{\cal E}_q \left( -; 0,{ib\over2} q^{(m+n)/2} \right)&
\left({ib\over2}\right)^n\cr
\times H_n & (x, q^{m/2} \vert  q)\ .}\eqno(43)
$$
The constructive derivation of this 
identity illustrates well the usefulness of
algebraic techniques for obtaining and interpreting properties of continuous
$q$-polynomials.  We plan to pursue investigations in this direction.

\vskip 2cm

\centerline{\bf Acknowledgements}
\smallskip
One of us (L.V.) is thankful to INFN, Sezione di Trieste, for hospitality and
support.

\vfill\eject

\centerline{\bf References}
\medskip

\item{1.}  Gasper, G. and Rahman, M., {\it Basic Hypergeometric Series},
(Cambridge University Press, Cambridge, 1990)
\smallskip
\item{2.} Koekoek, R. and Swarttouw, R.~F., The Askey-scheme of
hypergeometric orthogonal polynomials and its $q$-analogue, Report 94--05,
Delft University of Technology (1994).
\smallskip
\item{3.} Floreanini, R. and Vinet, L., A model for the continuous
$q$-ultraspherical polynomials, CRM--2233, Universit\'e de Montr\'eal (1995)
\smallskip
\item{4.} Kalnins, E.~G. and Miller, W., Symmetry techniques for $q$-series:
Askey-Wilson polynomials, Rocky Mountain J. Math. {\bf 19} (1989), 223--230
\smallskip
\item{5.} Floreanini, R. and Vinet, L., $q$-Orthogonal polynomials
and the oscillator quantum group, Lett. Math. Phys. {\bf 22} (1991), 45--54
\smallskip
\item{6.} Ismail, M.~E.~H. and Zhang, R., Diagonalization of certain
integral operators, Adv. Math. {\bf 109} (1994), 1--33.
\smallskip
\item{7.} Al-Salam, W., A characterization of the Rogers $q$-Hermite
polynomials, University of Alberta, preprint, 1994

\bye